\documentclass{amsart}

\usepackage{amsmath}

\usepackage{subfig}
\usepackage{mathrsfs}
\usepackage{amscd}
\usepackage{amscd, amssymb, amsmath, amsthm, graphics}
\usepackage{amsmath,amsfonts,amsthm,amssymb}
\usepackage{latexsym,amsmath}
\usepackage{graphicx,psfrag}
\usepackage[all]{xy}
\usepackage{latexsym}
\usepackage{enumerate}
\usepackage{verbatim}
\usepackage{pb-diagram}
\usepackage{accents}

\newtheorem{theorem}{Theorem}[section]
\newtheorem{lemma}[theorem]{Lemma}
\newtheorem{prop}[theorem]{Proposition}
\newtheorem{cor}[theorem]{Corollary}

\theoremstyle{definition}
\newtheorem{definition}[theorem]{Definition}

\theoremstyle{remark}
\newtheorem{remark}[theorem]{Remark}

\newtheorem*{claim}{Claim}

\numberwithin{equation}{section}

\newcommand{\AAA}{\mathbb{A}}
\newcommand{\CC}{\mathbb{C}}
\newcommand{\RR}{\mathbb{R}}

\newcommand{\TT}{\mathbb{T}}
\newcommand{\ZZ}{\mathbb{Z}}

\newcommand{\JJ}{\mathcal{J}}

\newcommand{\OO}{\mathcal{O}}
\newcommand{\MM}{\mathcal{M}}
\newcommand{\mmm}{\mathfrak{m}}

\def\R{\mathbb{R}}
\def\Z{\mathbb{Z}}
\def\T{\mathbb{T}}
\def\CP{\mathbb{CP}}
\def\RP{\mathbb{RP}}

\def\cM{\mathcal{M}}
\def\cN{\mathcal{N}}

\def\mf{\mathfrak}

\def\w{\omega}

\def\xkm2{\overline{X}_{k-2}}
\def\wt{\widetilde}

\newcommand{\bthm}{\begin{theorem}}
\newcommand{\ethm}{\end{theorem}}
\newcommand{\brmk}{\begin{remark}}
\newcommand{\ermk}{\end{remark}}

\newcommand{\thmone}{\begin{theorem}}
\newcommand{\thmtwo}{\end{theorem}}
\newcommand{\lemmaone}{\begin{lemma}}
\newcommand{\lemmatwo}{\end{lemma}}
\newcommand{\pfone}{\begin{proof}}
\newcommand{\pftwo}{\end{proof}}
\newcommand{\defone}{\begin{definition}}
\newcommand{\deftwo}{\end{definition}}
\newcommand{\corone}{\begin{cor}}
\newcommand{\cortwo}{\end{cor}}
\newcommand{\cone}{\begin{claim}}
\newcommand{\ctwo}{\end{claim}}
\newcommand{\propone}{\begin{prop}}
\newcommand{\proptwo}{\end{prop}}
\newcommand{\eqone}{\begin{equation}}
\newcommand{\eqtwo}{\end{equation}}
\newcommand{\rmkone}{\begin{remark}}
\newcommand{\rmktwo}{\end{remark}}
\newcommand{\enone}{\begin{enumerate}}
\newcommand{\entwo}{\end{enumerate}}
\newcommand{\itone}{\begin{itemize}}
\newcommand{\ittwo}{\end{itemize}}

\newcommand{\onehalf}{\left(\begin{array}{cc}}
\newcommand{\theother}{\end{array}\right)}

\newcommand{\oneeq}{\begin{equation}}
\newcommand{\twoeq}{\end{equation}}
\newcommand{\beq}{\begin{equation}}
\newcommand{\eeq}{\end{equation}}
\newcommand{\nono}{\noindent}

\newcommand{\tadj}{\JJ^{\epsilon}_{tadj}}
\newcommand{\MMc}{\MM_1^{crit}(2[D_4]; e_4,\widetilde X_0, L;\widetilde J)}
\newcommand{\MMb}{\MM_{k+1,l}^{\text{main}}(\beta, L)}
\newcommand{\MMm}{\MM_{k+1, l}^{\text{main}}}

\newcommand{\wkbc}{\widehat\cM_{\text{weak}}(L)}
\newcommand{\PO}{\mf{PO}}

\begin{document}

\title{On an exotic Lagrangian torus in $\CC P^2$}         
\author{Weiwei Wu}

\maketitle

\noindent{\small {\small\bf Abstract:} We find a non-displaceable
Lagrangian torus fiber in a semi-toric system which is superheavy
with respect to certain symplectic quasi-state.  In particular, this
proves Lagrangian $\RR P^2$ is not a stem in $\CC P^2$, answering a
question of Entov and Polterovich.  The main technique we apply is the relation
between Lagrangian Floer cohomology and symplectic quasi-morphisms/states due to
Fukaya, Oh, Ohta and Ono.\\

\vspace{2mm}\baselineskip 10pt

\noindent{\small\bf MR(2000) Subject Classification:}  53D12;
53D05\\

\section{Introduction}\label{section:Motivation}

The primary goal of this paper is to understand a toric degeneration
model of $\CC P^2$.
Our $\CP^2$ model should be considered as a
$\Z_2$-equivariant version of the one used in \cite{FOOOdeg} for
$S^2\times S^2$. However, we take a slightly different point of
view, based on symplectic cuts on cotangent bundles of manifolds
with periodic geodesics.  This degeneration gives a genuine torus action on an open part
of $\CP^2$, which results in an interesting family of Lagrangian
torus.  In particular, such degenerated torus action still gives a
moment polytope.  For $S^2\times S^2$, the polytope reads
$P_{S^2\times S^2}=\{(x_1,x_2)\in\RR^2:x_1,x_2\geq0,
x_1+2x_2\leq2\}\subset\RR^2$ (Figure $1$, see \cite{FOOOdeg}).  In
the case of $\CP^2$, one similarly have a toric action on an open
set, which gives a moment polytope as in Figure $2$, and can be
described as $P_{\CC P^2}=\{(x_1,x_2)\in\RR^2:x_1,x_2\geq0,
x_1+4x_2\leq4\}$.

\begin{figure}[h]
\begin{minipage}[t]{0.5\linewidth}
\centering
\includegraphics[width=1.5in]{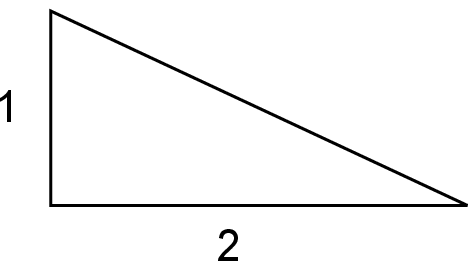}
\caption{$S^2\times S^2$} \label{fig:side:a}
\end{minipage}%
\begin{minipage}[t]{0.5\linewidth}
\centering
\includegraphics[width=1.5in]{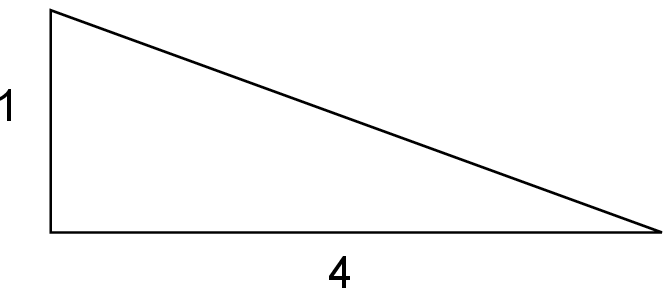}
\caption{$\CC P^2$} \label{fig:side:b}
\end{minipage}
\end{figure}

In \cite{FOOOdeg}, Fukaya-Oh-Ohta-Ono considered the Floer theory of
smooth fibers in the toric degeneration model of $S^2\times S^2$,
proving by bulk deformation that there are uncountably many
non-displaceable fibers in Figure $1$.  In view of
Albers-Frauenfelder's result \cite{AF}, we may interpret this result
as that, only non-displaceable torus fibers below the ``monotone
level" in the semi-toric system survives the symplectic cut along a
level set of $T^*S^2$.  This implies that the anti-diagonal of
$S^2\times S^2$ is not a stem (see \ref{section:review} for
the definition of a stem), answering a question raised by Entov and
Polterovich in \cite{EP2}. This was also proved independently by
several other authors \cite{ElP,CS}. In \cite{ElP} it was mentioned
that Wehrheim also has an unpublished note on this problem.

From Fukaya-Oh-Ohta-Ono's calculation on $S^2\times S^2$, we expect
the above similar picture of $\CC P^2$ also contains uncountably
many non-displaceable fibers.  This would correspond to an easy
adaption of Albers-Frauenfelder's result to $T^*\RR P^2$, by
considering the $\ZZ_2$-involution induced by antipodal map on
$S^2$.

In this paper we find one smooth non-displaceable monotone torus
fiber in the moment polytope described above, and prove that it is
superheavy with respect to some symplectic quasi-state.  In
particular, we proved:

\thmone\label{thm:main} There is a smooth monotone Lagrangian torus
fiber in Figure $2$, which is superheavy with respect to certain
 symplectic quasi-state.  In particular, it is stably
non-displaceable. \thmtwo

The limitation to such a monotone fiber is due to the difficulty of
using $\RR P^2$ as a bulk to deform our Floer cohomology as in \cite{FOOOdeg},
where the bulk is chosen to be an embedded $S^2$.  Since $\RP^2$
only obtains a non-trivial $\ZZ_2$-class, but our
calculation shows it is essential that we do not use coefficient ring of
characteristic $2$ (see Section \ref{s:completion}), the bulk of
$\RP^2$ does not improve our situation in a straightforward way. The author
does not know whether this is only technical. Nonetheless, our
computation suffices to show the following:

\corone $\RR P^2\subset \CC P^2$ is not a stem. \cortwo

\nono This answers the question of Entov-Polterovich (\cite{EP3},
Question 9.2) regarding the case of $\CC P^2$.

\rmkone It seems possible that our exotic monotone Lagrangian is in
fact the Chekanov torus in $\CC P^2$.  In particular, they both
bound $4$ families of disks of Maslov index $2$.  It would be nice
if one could identify the two geometrically, provided the guess is true.
 Nonetheless, even if such an identification holds, our calculation still
 gives new information: we would have an identification of Chekanov torus
 with a semi-toric fiber and showed the superheaviness of it.

 After this work was accomplished, Renato Vianna also detected certain exotic Lagrangian
 tori in $\CP^2$ \cite{Vianna}.  His tori can be discerned from the one presented here by the families
 of holomorphic disks with Maslov number $2$ they bound.  It would be interesting to compare
 the two different objects.\rmktwo

\nono\textit{Acknowledgement:}  The author would like to thank
Yong-Geun Oh and Leonid Polterovich for their interests in this
work, and many inspiring discussions, which deepens the author's
understandings to this question.  The author is indebted to Kenji
Fukaya for many interesting comments, as well as pointing out an
error in an early draft.  The author is also grateful to Chris Wendl
for patiently explaining many details of the automatic
transversality, to Garrett Alston for drawing the author's attention
to \cite{FOOOdeg}, to Matthew Strom Borman for explaining details of
symplectic quasi-states, and to Jie Zhao, Ke Zhu for many useful
discussions.  Part of this work was completed when the author was a
graduate student in University of Minnesota under the supervision of
Tian-Jun Li. The author is supported by FRG 0244663.

\section{Preliminaries}\label{section:review}

The current section summarizes part of the Lagrangian Floer theory
developed by Fukaya-Oh-Ohta-Ono \cite{FOOO,FOOOtoricI,FOOOtoricII}
etc, as well as the theory of symplectic quasi-states developed by
Entov and Polterovich in a series of their works \cite{EP1,EP2,EP3}
etc. The aim of this section is to recall basic notions and main
framework results in these theories for our applications, as well as
for the convenience of readers. Therefore, our scope is rather
restricted and will not provide a thorough account to the whole
theory. For details and proofs one is refered to the above-mentioned
works. Much of our discussions on Lagrangian Floer theory follow the
lines of \cite{FOOOdeg}.

\subsection{Lagrangian Floer theory via potential function}\label{subsection:potential function}

Let $(M,\omega)$ be a smooth symplectic manifold and $L\subset M$ a
relatively spin Lagrangian.  This means the second Stiefel-Whitney
class $w_2(L)$ is in the image of the restriction map
$H^2(M,\ZZ_2)\rightarrow H^2(L,\ZZ_2)$. We first describe the moduli
spaces under consideration.  Let $J\in \JJ_\omega$, the space of
compatible almost complex structures, and $\beta\in H_2(M,L;\ZZ)$.
We denote by $\MM^{\text{main}}_{k+1,l}(\beta;M,L;J)$ as the space
of $J$-holomorphic bordered stable maps in class $\beta$ with $k+1$
boundary marked points and $l$ interior marked points. Here, we
require the boundary marked points to be ordered
counter-clockwisely. When no confusions is likely to occur, we will
suppress $M$, $L$ or $J$.

One of the fundamental results in \cite{FOOO} shows that, one has a Kuranishi structure on $\MMb$, so that
the evaluation maps at the $i^{\text{th}}$ boundary marked point ($j^{\text{th}}$ interior marked point, respectively)

$$ev_i:\MMb\rightarrow L,$$

\nono and

$$ev_j^+:\MMb\rightarrow M$$

\nono are weakly submersive (see \cite{FOOO} for the definition of
weakly submersive Kuranishi maps). For given smooth singular
simplices $(f_i:P_i\rightarrow L)$ of $L$ and $(g_j:Q_j\rightarrow
M)$ of $M$, one can also define the fiber product in the sense of
Kuranishi structure:

$$\MM_{k+1,l}^{\text{main}}(\beta;L;\overrightarrow{Q},\overrightarrow{P}):=\MM_{k+1, l}^
{\text{main}}(\beta;L)_{(ev_1^+,\dots,ev_l^+,ev_1,\dots,ev_k)\times(g_1,\dots,g_l,f_1,\dots,f_k)}(\Pi_{j=1}^l
Q_j\times \Pi_{i=1}^k P_i).$$

The virtual fundamental chain associated to this moduli space,
$$ev_0:\MM_{k+1,l}^{\text{main}}(\beta;L;\overrightarrow{Q},\overrightarrow{P})\rightarrow L$$
as a singular chain, is defined in \cite{FOOO} via techniques of virtual perturbations.



We consider the universal Novikov rings:

$$\Lambda=\{\sum a_iT^{\lambda_i}| a_i\in\CC, \lambda_i\in\RR, \lambda_i\leq \lambda_{i+1}, \lim_{i\rightarrow\infty}\lambda_i=\infty\},$$
$$\Lambda_{0}=\{\sum a_iT^{\lambda_i}| a_i\in\CC, \lambda_i\in\RR_{\geq0}, \lambda_i\leq \lambda_{i+1}, \lim_{i\rightarrow\infty}\lambda_i=\infty\}.$$

Here $T$ is a formal variable.  Consider a valuation which assigns
$\sigma_T(\sum a_iT^{\lambda_i})=\lambda_1$ if not all $a_i=0$, and
let $\sigma_T(0)=+\infty$.  This induces a $\RR$-filtration on
$\Lambda$ and $\Lambda_0$ thus a non-Archimedian topology.  Note
that $\Lambda_0\subset\Lambda$, and $\Lambda_0$ has a maximal ideal
$\Lambda_+$ consisting of elements with $\lambda_i>0$ for all $i$.
The absence of $e$-variable will reduce the grading of Floer
cohomology groups to $\Z_2$, but this is irrelevant to our
applications.

The heart of Fukaya-Oh-Ohta-Ono's work is to define a filtered
$A_\infty$-structure on $C^*(L;\Lambda_0)$ for a Lagrangian
$L\subset M$ and define Floer cohomology by deformations of (weak)
bounding cochains.  However, we will not mention the explicit
constructions here, both because they are far beyond our scope, and
that there are plenty of comprehensive reference and surveys
available in the literature. Just for an incomplete list:
\cite{FOOO, FOOOtoricI, FOOOtoricII, FOOOsurvey}, etc.  Instead we
will adopt a most economic approach towards the applications in
mind, by recalling a package made available by the deep theory,
namely, the computations on Lagrangian Floer cohomology via
potential functions.

A \textit{potential function} $\frak{PO}^L$ is a $\Lambda_+$-valued
function defined on the set of \textit{weak bounding cochain} of
$L$, denoted as $\widehat\cM_{\text{weak}}(L)$.  For any
$b\in\wkbc$, one may associate a Floer cohomology group $HF^*(L,b)$
for the pair $(L,b)$.  We do not define the weak bounding cochains
in general, however, according to \cite[Remark A.2]{FOOOdeg},
$\wkbc$ can be identified with
$H^1(L;\Lambda_0)/H^1(L;2\pi\sqrt{-1}\Z)$ for any monotone
Lagrangian submanifolds with minimal Maslov number equal $2$.  Hence
in the rest of this paper, $\wkbc$ will refer to this particular
set, and the potential function can be written as:

\beq\label{eq:POdef} \mf{PO}^L:
H^1(L;\Lambda_0)/H^1(L;2\pi\sqrt{-1}\Z)\rightarrow \Lambda_+.\eeq

With the monotonicity assumption above, one may compute $\mf{PO}$
 explicitly as in \cite[Theorem A.1, A.2]{FOOOdeg}.  Choose a basis $\{ e_i\}_{i=1}^n$ for
$H^1(L;\Z)$ and represent $b=\sum_{i=1}^n x_ie_i $ for $b\in
H^1(L;\Lambda_0)$ and $x_i\in\Lambda_0$.  Then the potential
function is written as:

\beq\label{eq:POtemp}
\PO^L(b)=\sum_{\mu{\beta}=2}ev_0^*([\cM_1(L;J,\beta)])T^{\omega(\beta)/2\pi}exp(b(\partial\beta)).\eeq

Here $[\partial\beta]\in H_1(L;\Z)$, hence
$b(\partial\beta)\in\Lambda_0$.  Writing in coordinates, the
potential function can be regarded as a function from
$(\Lambda_0/2\pi\sqrt{-1}\Z)^n$ to $\Lambda_+$.  A change of
coordinate $y_i=e^{x_i}$ transforms the function in
\eqref{eq:POtemp} into the form commonly used in the literature:

\begin{equation}\label{eq:PO}
\begin{aligned}
\PO^L:&\hskip2mm(\Lambda_0\backslash \Lambda_+)^n&\rightarrow &\hskip3mm\Lambda_+\\
&(y_1,\dots,y_n)&\mapsto&\sum_{\mu(\beta)=2}ev_0^*([\cM_1(L;J,\beta)])T^{\omega(\beta)/2\pi}y_i^{l_i},
\end{aligned}
\end{equation}

where $\partial\beta=\sum_{i=1}^n l_ie_i^*$, for $\{e_i^*\}_{i=1}^n$
a dual basis in $H_1(L;\Z)$.  The following result manifests the
importance of potential functions:

\bthm[\cite{FOOOdeg}, Theorem 2.3]\label{thm:critical points} Let
$L$ be a Lagrangian torus in $(M,\w)$.  Suppose that
$H^1(L;\Lambda_0)/H^1(L;2\pi\sqrt{-1}\Z)\subset \wkbc$ and $b\in
H^1(L;\Lambda_0)$ is a critical point of the potential function
$\PO^L$ of $L$.  Then we have:

$$HF^*(L,b)\cong H^*(L;\Lambda_0).$$

\nono In particular, $L$ is non-displaceable. \ethm

\brmk Results state in this section holds valid for
Lagrangians satisfying Condition 6.1 in \cite{FOOOdeg}, that is, if
any non-empty moduli space of holomorphic disks with Maslov index
$\geq2$.  For monotone Lagrangians, there is also an alternative
approach developed by Biran and Cornea \cite{BC09} via pearl
complexes, which was used in Vianna's work \cite{Vianna} and many
others.

\ermk

\subsection{Symplectic quasi-states and Lagrangian Floer theory}

In this section we briefly review the theory of symplectic
quasi-states developed by Entov and Polterovich. A
\textit{symplectic quasi-state} is a functional
$\zeta:C^\infty(M)\rightarrow\RR$ satisfying the following axioms
for $H,K\in C^\infty(M)$ and $\lambda\in\RR$:

\enone[(i)]
\item(Normalization) $\zeta(1)=1$;
\item(Monotonicity) If $H\leq K$, then $\zeta(H)\leq \zeta(K)$;
\item(Quasi-linearity) If $\{H,K\}=0$, then $\zeta(H+\lambda K)=\zeta(H)+\lambda\zeta(K)$;
\item(Vanishing) If supp$(H)$ is displaceable, then $\zeta(H)=0$;
\item(Symplectic invariance) $\zeta(H)=\zeta(H\circ f)$ for $f\in Symp_0(M)$
\entwo

Given a symplectic quasi-state $\zeta$ and a subset $S\subset M$,
$S$ is called $\zeta$-\textit{heavy} if:
$$\zeta(F)\geq \inf_{x\in S}F(x), \hskip 2mm\forall F\in C^\infty(M);$$

\nono and $\zeta$-\textit{superheavy} if
$$\zeta(F)\leq \sup_{x\in S}F(x), \hskip 2mm\forall F\in C^\infty(M).$$

One of the basic properties of these subsets proved in \cite{EP2} is
that, a $\zeta$-superheavy subset is always $\zeta$-heavy, and a
$\zeta$-heavy set is stably non-displaceable (this is a notion strictly stronger
than non-displaceability).  Let $V\subset C^\infty(M)$ be a finite
dimensional linear subspace spanned by pairwisely Poisson-commuting
functions. Let $\Psi:M\rightarrow V^*$ be the moment map defined by
$\langle\Psi(x),F\rangle=F(x)$ for $F\in V$.  A non-empty fiber of
this moment map is called a \textit{stem}, if the rest of the fibers
are all displaceable.  It was essentially proved in Theorem 1.6 of
\cite{EP2} the following:

\thmone[\cite{EP2}] A stem is a superheavy subset with respect to
arbitrary symplectic quasi-states. \thmtwo

In general, the existence of symplectic quasi-states is already an
intriguing question.  \cite{EPsemi} showed that, given a direct sum
decomposition of $QH^*_{2n}(M)=\mathbb{F}\oplus QH'$, where
$\mathbb{F}$ is a field, then one may associate a symplectic
quasi-state $\zeta_e$ to the unit element $e\in\mathbb{F}$.

The relations between symplectic quasi-states and Lagrangian Floer
theory are established by the $i$-\textit{operator} (sometimes also
referred to as the open-closed string maps or the Albers map in the
literature).  The version of $i$-operator we need involves a
deformation by the weak bounding cochains, thus denoted as:

$$i^*_b:QH^*(M)\rightarrow HF(L,b).$$

The concrete definition of $i^*_b$ was given in \cite{FOOO}, and we
refer interested readers there for details (see also \cite{BC09} for
a similar operator in the context of pearl complexes).  The key
property of $i^*_b$ we need is that it sends the unit of $QH^*(M)$
to that of $HF(L,b)$.  This fact was shown in 7.4.2-7.4.6 in
\cite{FOOO}, which passes to the so-called \textit{canonical model}
of $C^*(L;\Lambda_0)$ and involved deep algebraic techniques in
filtered $A_\infty$ algebras, therefore is beyond the scope of the
present paper.

With this understood, the key results our proof will rely on reads
as follows:

\thmone[Theorem 18.8, \cite{FOOOspectral}]\label{thm:FOOO spectral}
Let $L$ be a relatively spin Lagrangian submanifold of $M$,
$b\in\widehat\MM_{\text{weak}}(L)$ be a weak bounding cochain.
$e\in QH^*(M;\Lambda)$.

\enone[(1)]
\item If $e\cup e=e$ and $i_{b}^{*}(e)\neq0$, then $L$ is
$\zeta_e$-heavy.
\item If $QH^*(M;\Lambda)=\Lambda\oplus Q$ is a direct factor
decomposition as a ring, and $e$ comes from a unit of the factor
$\Lambda$ which satisfies $i^*_{b}(e)\neq0$, then $L$ is
$\zeta_e$-superheavy. \entwo

\thmtwo

\corone\label{cor:main lemma} Suppose $
QH^*(M;\Lambda)=\bigoplus_{i=1}^n\Lambda e_i$ as a ring, for $e_i\in
QH^*(M;\Lambda)$ being a series of idempotents (in particular $QH^*$
is semi-simple).  If $HF^*((L,b);\Lambda)\neq0$, then $L$ is
superheavy for certain symplectic quasi-state $\zeta_{e_k}$, $1\leq
k\leq n$.

\cortwo

\pfone This is implicit from the proof of Theorem 23.4,
\cite{FOOOspectral}.  Since
$i^*_{qm,b}$ sends the unit to the unit, at least one of
the idempotents $e_k$ has non-vanishing image.  From Theorem
\ref{thm:FOOO spectral}, $L$ is superheavy. \pftwo



Combining Corollary \ref{cor:main lemma}, Theorem \ref{thm:critical
points} and \eqref{eq:PO}, provided we have a semi-simple quantum
cohomology ring for the ambient manifold $M$, to show a monotone
Lagrangian torus is superheavy with respect to certain symplectic
quasi-state, it suffices to compute the contribution of each moduli
space of holomorphic disks of Maslov index $2$, and find the critical points
for the potential function, which will be the
topic of subsequent sections.

\section{A semi-toric system of $\CC P^2$}\label{section:semitoric
system}

\subsection{Description of the system}\label{subsection:description of CP2}

We recall the semi-toric system of $\CC P^2$ particularly suitable
for our problem.  We first briefly recall the semi-toric model for
$S^2\times S^2$ following the idea of \cite{EP1} and \cite{Seidel}.
Write $S^2\times S^2$ as

$$\{x_1^2+y_1^2+z_1^2=1\}\times \{x_2^2+y_2^2+z_2^2=1\}\subset \RR^3\times \RR^3.$$
\nono Let
$$\widetilde{F}(x_1,y_1,z_1;x_2,y_2,z_2)=z_1+z_2,$$
$$\widetilde{G}(x_1,y_1,z_1;x_2,y_2,z_2)=\sqrt{(x_1+x_2)^2+(y_1+y_2)^2+(z_1+z_2)^2},$$
\nono then
$$\Phi_{S^2\times S^2}=(F,G):=(\widetilde{F}+\widetilde{G},2-\widetilde{G}):S^2\times S^2\rightarrow \RR^2$$

\nono defines a Hamiltonian system. $\widetilde G$ is not integrable
when it equals $0$, that is, at the anti-diagonal $\bar\Delta$. This
Hamiltonian system gives a moment polytope as in Figure $1$ up to a
rescale of the symplectic form, with a singularity at $(0,1)$
representing a Lagrangian sphere.  In classical terms, this is in
fact a moment polytope for $S^2\times S^2\backslash\bar\Delta$,
where any tubular neighborhood $N(\bar\Delta)$ of $\bar\Delta$ will be
mapped into to a neighborhood of $(1,0)$.

Another useful point of view is to consider $S^2$ equipped with the
standard round metric, which induces a metric on its cotangent
bundle. $S^2\times S^2$ is obtained from $T^*S^2$ by a symplectic
cut at the hypersurface
$$M_1=\{p\in T^*S^2: |p|=1\}.$$

The circle action on this hypersurface is exactly the unit-speed
geodesic flow we use for cutting.  See \cite{Lerman} for details of
the construction of symplectic cuts. In this perspective, we may
describe the $\mathbb{T}^2$-action induced by $\Phi$ in a geometric
way. Consider the rotation of $S^2$ along an axis. The cotagent map
of this rotation generates the circle action $\tau_{\widetilde{F}}$
on the whole $T^*S^2$. Another circle action $\tau_{\widetilde{G}}$
is generated by the unit geodesic flow on the complement of the zero
section (we already used it for symplectic cut above). Both
$\tau_{\widetilde{F}}$ and $\tau_{\widetilde{G}}$ descend under the
symplectic cut and commute, thus induces a genuine $\TT^2$-action on
$S^2\times S^2\backslash \bar\Delta$.\\

We proceed to the case for $\CC P^2$.
Consider the $\ZZ_2$-action on $T^*S^2$ induced by the antipodal map
on the zero section.  
It is readily seen that, the symplectic cut at the level set $M_1$ is also
$\ZZ_2$-equivariant, so we may well quotient out this $\Z_2$ action first and
then perform the symplectic cut. This is equivalent to performing symplectic cut
on $T^*\RP^2$, which results in a symplectic $\CP^2$.  In summary
 we have the following commutative
diagram, which is equivariant with respect to the action of
$\tau_{\tilde F}$ and $\tau_{\tilde G}$:

$$\xymatrix{
&T_1^*S^2\ar@{^{(}->}[r]\ar^{\pi}@{->}[d] &S^2\times S^2\ar^{\iota}@{->}[d]\\
&T^*_1\RR P^2\ar@{^{(}->}[r] &\CC P^2}$$

Here $\pi$ is the $2$-to-$1$ cover over $T^*\RR P^2$ and $\iota$ is the
standard two-fold branched cover from $S^2\times S^2$ to $\CC P^2$,
branching along the diagonal. Notice now both $\tau_{\tilde F}$ and
$\tau_{\tilde G}$ are $\Z_2$-equivariant under the deck
transformation, therefore, the above $2$-fold cover induces two
commuting circle action on $\CP^2\backslash \RP^2$.  However, since
the $\Z_2$-action halves the length of each geodesic, to get a
time-1 periodic flow, the Hamiltonian function generating the
circle action descended from $\tau_{\tilde G}$ should be the
descendant of $\frac{1}{2}\tilde G$ on $S^2\times S^2$.  The end
result after approapriate reparametrizations is a toric model for
$\CP^2\backslash Q$ with moment polytope as in Figure $2$. Note from
the reasoning regarding $\tilde G$, after the reparametrization the
line area of $\CC P^2$ is $2$ (if the line class area had been $1$, the
sizes in Figure $2$ would have been $\frac{1}{2}$ by $2$). Similar
to the case of $S^2\times S^2$, $(1,0)$ indeed represents the
standard Lagrangian $\RP^2\subset \CP^2$. We will denote this
semi-toric moment map as $\Phi_{\CC P^2}$.

\subsection{Symplectic cutting $\CC P^2$}\label{subsection:cutting
CP2}

The main ingredient of our proof, following an idea of the arxiv
version of \cite{FOOOdeg}, is to split $\CC P^2$ into two pieces and
glue the holomorphic curves. The splitting we use is described as
follows. We continue to regard $\CC P^2$ as a result of cutting
along $M_1$ in $T^*\RR P^2$. Consider
$M_\epsilon=\{|p|=\epsilon\}\subset T^*\RR P^2\hookrightarrow \CC
P^2$. A further symplectic cut along $M_\epsilon$ results in two
pieces, and we examine this cutting in slightly more detailed.

Let $X_0$, $X_1$ be the two components of $\CC P^2\backslash
M_\epsilon$, where $X_1$ contains the original $\RR P^2$. Their
closure, denoted $X'_0$ and $X'_1$, respectively, has a boundary
being the lens space $L(4,1)$ equipped with the standard contact
form (the one coming from $S^3$ quotiented by a $\ZZ_4$-action), and
therefore a local $S^1$-action of the neighborhood.

As is constructed in \cite{Lerman}, by quotienting such an action on
$\partial X_1$ and gluing back to $X'_1$, one completes the
symplectic cutting and this operation results in $X''_1:=(\CC
P^2,2\epsilon\omega_0)$. Denote $H\in H_2(\CP^2,\Z)$ as the homology class of
a line, then $X''_1\backslash X_1$ is an embedded
symplectic divisor in $X''_1$ of class $2H$, which we called the
\textit{cut locus} or \textit{cut divisor}.  The same procedure on the other piece $X_0$
leads to a minimal symplectic $4$-manifold (see for example Lemma
$1.1$ in \cite{Seppi}), along with a symplectic sphere of
self-intersection $(+4)$ inherited from the quadric
$Q:=\{x^2+y^2+z^2=0\}$ in the original $\CC P^2$. Moreover, it
contains a symplectic sphere of self-intersection $(-4)$ as cut
locus, from which we see $X''_0$ is indeed the symplectic fourth
Hirzebruch surface $F_4$ by McDuff's famous classification of
rational and ruled manifolds \cite{Mcduff}.

We also want to examine such a cut process from the other side of
$M_1$. Biran's decomposition theorem for $\CC P^2$ (\cite{Biran
barrier}) implies that $\CC P^2\backslash \RR P^2$ is indeed a
symplectic disk bundle $\OO(4)$ over a sphere, where the zero
section has symplectic area $4$, and the symplectic form is given by
$\pi^*\omega_{\Sigma}+d(\bar r^2\alpha)$.  Here $\pi$ is the
projection to the zero section, $\omega_{\Sigma}$ a standard
symplectic form on the sphere up to a rescale, $\bar r$ the radial
coordinate of the fiber and $\alpha$ a connection form of the circle
bundle associated to $\OO(4)$. Then the fiber class has at most
symplectic area $1$, and the total space can be identified
symplectically with $\CC P^2\backslash \RR P^2$ with the standard
sympletic form.

In this case $X'_0$ is identified with $\{|\bar
r|\leq1-\epsilon\}\subset \OO(4)$ and the geodesic flow in $T^*\RR
P^2$ is identified with the action of the one obtained by
multiplying $e^{i\theta}$ in each fiber.  Therefore, one may perform
a symplectic cut along $|\bar r|=1-\epsilon$ for $1\gg\epsilon>0$,
the resulting manifold is again the symplectic $F_4$ as above, where the form is
compatible with the standard (integrable) complex structure obtained
as $P(\OO\oplus\OO(4))$. To summarize, we have the following (see
also Figure $3$):

\lemmaone Consider $\CC P^2$ as a consequence of symplectic cut
along the contact type hypersurface $\{|p|=1\}\subset T^*\RR P^2$.
Then a further symplectic cut along $\{|p|=\epsilon\}$ results in a
$\CC P^2$ with rescaled symplectic form, as well as a symplectic
fourth  Hirzebruch surface whose zero section has symplectic area
equal $2$.  Moreover, the symplectic $\CC P^2$ comes naturally with
a symplectic quadric as the cut locus, and $F_4$ with a
$(-4)$-sphere as cut locus.\lemmatwo

\rmkone Discussions above seem to be well-known.  For a dual
perspective via symplectic fiber sum, one is referred to for example
\cite{Seppi}.  In particular, the above cutting can be seen as a
reverse procedure of symplectic rational blow-down of the
$(-4)$-sphere in the symplectic $F_4$. \rmktwo

\begin{figure}[h]\label{fig:cutting along M}
\centering
\includegraphics[width=4.5in, height=2in]{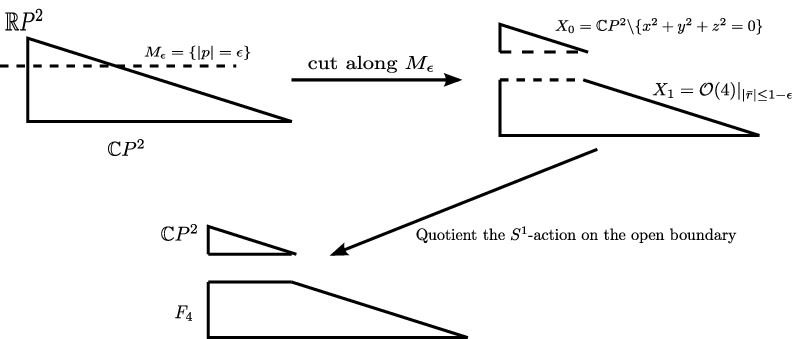}
\caption{Cutting along $M_\epsilon$}
\end{figure}

\subsection{Second homology classes of $\CC P^2$ with boundary on a semi-toric
fiber}\label{subsection:description of classes}

From Section \ref{subsection:description of CP2}, we have obtained a desired family of
Lagrangian torus as semi-toric fibers in $\CP^2$.  From now on $L$
will denote one of the semi-toric fibers parametrized by
$\R^2$-coordinates in Figure 2.  Our next task is to understand
$H_2(\CP^2,L)$. 

From the usual long exact sequence for relative homology, one easily
sees that $H_2(\CC P^2,L;\ZZ)$ has rank $3$.  Again split $\CP^2$ along $M_\epsilon$
into a copy of $F_4$ and $(\CP^2,\epsilon\w_{std})$
as in Section \ref{subsection:cutting CP2}, while keeping $L\subset F_4$ by choosing
$\epsilon$ small enough.

From the classification theorem of homology classes in \cite{CO},
one obtains eight homology classes
of interests, marked as $[D_i]$ and $[e_i]$, $i=1,2,3,4$ as Figure
$4$ below.  In the figure, $e_i$ are the $\T^2$-equivariant divisors
and, as relative cycles, $D_i$ denote the image of $J_0$-holomorphic
disks which intersects $e_j$ exactly $\delta_{ij}$ times counting
multiplicity. For ease of drawing we did not draw $D_3$
perpendicular to $e_3$, but it is understood in the way of how Cho
and Oh described in \cite{CO}.

\begin{figure}[h]
\centering
\includegraphics[width=3.7in, height=1.2in]{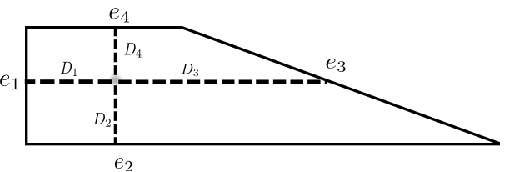}
\caption{$H_2(F_4,L;\ZZ)$}
\end{figure}

 Out of these eight classes, one has a basis of
$H_2(F_4,L;\ZZ)$ consisting of $[e_1]$, $[e_2]$, $[D_1]$ and
$[D_2]$.  Other classes have relations

\begin{align} \label{eq:rel1}    &[e_3]=[e_1], \hskip 2mm [e_4]=[e_2]-4[e_1],\\
                    \label{eq:rel2}   &[D_1]+[D_3]+4[D_2]=[e_2], \hskip 2mm [D_2]+[D_4]=[e_1].\end{align}

On way of checking these relations is to use Poincare pairings and
gluing chains with opposite boundaries on $L$.  Notice that there is
a natural homomorphism by restriction to the Borel-Moore homology of $F_4\backslash e_4$:

$$\iota: H_2(F_4,L;\ZZ)\rightarrow H_2^{BM}(F_4\backslash e_4, L;\ZZ).$$

\nono$\iota$ is surjective with kernel $[e_4]$, so the classes in
$H_2^{BM}(F_4\backslash e_4, L;\ZZ)$ can still be represented by
$e_i$, $i=1,2,3$ and $D_j$, $j=1,2,3,4$ appropriately punctured with
the same relations as in (\ref{eq:rel1}), (\ref{eq:rel2}). These
facts can be easily seen from the duality between the Borel-Moore
homology and the usual cohomology.

On the other side of the cutting, which is $\CC P^2\backslash Q$,
where $Q=\{x^2+y^2+z^2=0\}$ being the standard quadric, the second
Borel-Moore homology contains only a $2$-torsion.
We will only consider Borel-Moore cycles with
asymptotics equal a union of certain Reeb orbits of
$\partial^\infty(\CC P^2\backslash Q)=L(4,1)$.  Regard Borel-Moore cycles with $2k$ punctures
(number of Reeb orbits)
at infinity as equivalent, and denote such equivalence classes by $kH'$.
Note that $k$ already contains information of the homology classes: cycles in $k_1H'$ and
$k_2H'$ represent the same Borel-Moore classes in
$H_2(\CP^2\backslash Q,\partial^\infty(\CP^2\backslash Q))$ if and
only if $k_1-k_2\equiv0$ mod$2$, but the relative chern number will
depend on the actual equivalence classes instead of solely the
Borel-Moore classes.  See Section \ref{subsection:CZ indices and
c1}.
\\

\nono{\bf Relations between classes in $X_0=F_4\backslash e_4$ and
$X_1=\CC P^2\backslash Q$:}
 To describe the relations between classes in the two pieces, we first fix a basis
of $H_2^{BM}(F_4\backslash e_4, L;\ZZ)$ consisting of
$\{\iota[e_2], \iota[D_1], \iota[D_2], \iota[D_4]\}$.  When no
possible confusion occurs, we will simply suppress $\iota$ by abuse of notation. Notice
that cycles in $kH'$ in $X_1$ has $2k$ punctures counting
multiplicity, which matches with cycles with coefficient $2k$ in the
$D_4$-component in $X_0$.  Of particular interests, by matching a
cycle $C_{H'}\subset X_1$ in class $H'$ with a $2$-cycle of class
$2[D_4]$ with correct asymptotics, one obtains a relative cycle in
$\CC P^2$ with boundary on $L$.  The class in $H_2(\CC P^2, L;\ZZ)$
represented by such a cycle is denoted as $[D'_4]=2[D_4]\#[H']$.

To understand $[D_4']$ more explicitly, notice that
$\partial[D'_4]=2\partial[D_4]=-2\partial[D_2]\in H_1(L;\ZZ)$.
Therefore, one may match a cycle in class $[D'_4]$ with one in
$2[D_2]$ to obtain a closed cycle in $\CC P^2$.  Such a cycle
intersects $e_2$ positively twice counting multiplicities, and
therefore represents nothing but class $H\in H_2(\CC P^2;\ZZ)$. In
summary, we deduced that:

\eqone\label{eq:D4 and D2 consist H}
[H']\#2[D_4]\#2[D_2]=[D_4']\#2[D_2]=H\in H_2(\CC P^2;\ZZ)\eqtwo

Classes $[e_1]$ and $[e_3]$ does not extend naturally to closed
classes as in $H_2(\CC P^2)$.  However, as what we did to $[D_4]$,
twice of them caps cycles in $H'$ of $X_1$.  Therefore we also have:

\eqone\label{eq:e13 and H' consist H} [H']\#2[e_1]=[H']\#2[e_3]=H\in
H_2(\CC P^2;\ZZ)\eqtwo

These gluing relations will play an important role later.  It is
also readily seen that $\{H,[D_1],[D_2]\}$ forms a basis of $H_2(\CC
P^2, L;\ZZ)$, where $F_4\backslash e_4$ is (symplectically) embedded
to $\CC P^2$ in a canonical way, thus induces a natural inclusion of
Borel-Moore two cycles.

\subsection{Computation of the relative Chern numbers and Conley-Zehnder
indices}\label{subsection:CZ indices and c1}

We now compute the Maslov indices for $H_2(\CP^2,L)$ by
understanding the relative Chern classes and Conley-Zehnder indices
involved. A technical reason which makes our case slightly more
complicated than the case of a Lagrangian $S^2$ is that there is no
natural splitting of $T(T^*\RR P^2)$. This is caused by the
non-orientability of $\RR P^2$ (to compare the case of Lagrangian
$S^2$, see for example \cite{Hind, EvansT, LW12}).  However, we will use a
trivialization of the splitting surface $M_\epsilon$ which seems
even more natural and convenient in the (semi-)toric context.

As we already saw, there is an $S^1$-action on $\partial
X'_i=M_\epsilon$ for both $i=0,1$.  In the toric picture of $\CC
P^2\backslash \RR P^2$, such an $S^1$-action induces a vector field
on $M_\epsilon$ which is dual to $\frac{\partial}{\partial x_2}$ in
the moment polytope.  This action induces a natural trivialization
of the contact distribution over its own orbits.  We will call such
a trivialization $\Phi$ and use it to compute the Conley-Zehnder
indices and first Chern numbers.  For the definitions of these two
invariants one is referred to \cite{SFT}, or \cite{EvansT,Hind}

By definition, the Poincare return map with respect to such a
trivialization is always identity, therefore,

\eqone\label{eq:CZ vanishes}\mu^\Phi_{CZ}\equiv 0.\eqtwo

We will pursue the first Chern number for (Borel-Moore) classes
described in Section \ref{subsection:description of classes} in the
rest of the section.

We start with $X_0$.  As always we assume the Lagrangian torus fiber
is contained in this side.  Consider again the $\OO(4)$ disk bundle
as in Section \ref{subsection:cutting CP2}, from which we cut along
another hypersurface $M_{\epsilon/2}=\{r=1-\frac{\epsilon}{2}\}$ to
obtain a symplectic fourth Hirzebruch surface $\overline X_0$.  One
may also equip it a compatible toric complex structure.  The
anti-canonical divisor is defined by the equivariant divisors on the
boundary of the moment polytope, therefore, the anti-canonical line
bundle $\bigwedge^2 T\overline X_0$ admits an equivariant section
$\xi$ vanishing exactly on the boundary equivariant divisors with
order $1$.

Embed $X_0$ equivariantly into $\overline{X}_0$.  Take any cycle
$u:\Sigma\rightarrow X_0$ with boundary on a torus fiber $L$ and
asymptotics being Reeb orbits of $\partial X_0$.  It has boundary
Maslov index zero if we take the trivialization induced by the torus
action near $L$. Assume that $u$ intersects transversally with the
equivariant divisors. The pull-back $u^*\bigwedge^2 (TX_0,J)$ thus
comes naturally with a section $u^*\xi$ which vanishes at
$u(\Sigma)\cap \bigcup_{i=1}^4 e_i$ with order $\pm1$ depending on
the intersection form.  $u^*\xi$ is clearly equivariant with the
$S^1$-action on $\partial X_0$ and the torus boundary thus agreeing
with the trivialization there. This observation computes immediately
the following:

\eqone\label{eq:first cherns 1}
c_1^\Phi(D_1)=c_1^\Phi(D_2)=c_1^\Phi(D_3)=1, \hskip 2mm
c_1^\Phi(D_4)=0. \eqtwo

Notice also that the first chern number of $e_2$ is independent of
the choice of trivializations.  From (\ref{eq:rel1}) and
(\ref{eq:rel2}) we may compute the rest of the chern numbers
summarized as follows:

\eqone\label{eq:first cherns 2} c_1^\Phi(e_1)=c_1^\Phi(e_3)=1,
\hskip 2mm c_1^\Phi(e_2)=6. \eqtwo

For the relative Chern classes in $X_1$, we again focus on cycles
with asymptotics equal copies of $S^1$-orbits on $M_\epsilon$. Note
that, when counting multiplicity, there are always even number of
$S^1$-orbits since simple orbits represents a non-trivial element in
$\pi_1(T^*\RR P^2)$. The class $kH'$ has $2k$ asymptotics, which can
be capped by $2k[D_4]\#2k[D_2]$ to form a closed cycle in $\CC
P^2$ from \ref{eq:rel2}. Such a class intersects positively with
$e_2$ at $2k$ points, thus itself being the class $kH$ in $\CC P^2$.
From our computation in $X_0$, we see that

\eqone\label{eq:first cherns 3} c_1^\Phi(kH')=3k-2k=k>0.\eqtwo

\section{Classification of Maslov $2$ disks}

\subsection{A quick review on SFT and neck-stretching}\label{section:neck-stretching}

In this section we collect basic definitions and facts from
symplectic field theory, especially the part of neck-stretching,
mostly for readers' convenience and to fix notations.  For more
details, we refer interested readers to \cite{SFT},
\cite{compactness}, and other expositions such as \cite{EvansT,
Hind, LW12}.

Given a closed symplectic manifold $(M,\omega)$, we call
$(N,\alpha)$ a \textit{contact type hypersurface} if there is a
neighborhood $V$ of $N$, such that $V$ is diffeomorphic to
$(-\epsilon, \epsilon)\times N$, and $\partial_s$ is a Liouville
vector field in $V$, that is,
$\mathcal{L}_{\partial_s}\omega=\omega$. Here $s$ is the coordinate
of the first component of $U$. In this case,
$\alpha=i_{\partial_s}\omega$ is a contact form, of which the
contact distribution is denoted $\xi$, and the Reeb flow denoted by
$R$.

An almost complex structure $J\in \mathcal{J}_\omega$ is called
\textit{adjusted} if the following conditions hold in $U$:

\begin{enumerate}[(i)]
\item $J|_\xi=\widetilde{J}$ is independent of $s$;

\item $J(\partial_s)=R$.
\end{enumerate}
\vskip 2mm

We now consider a deformation of a given adjusted almost complex
structure $J$.  Let $V_t=[-t-\epsilon,t+\epsilon]$ and
$\beta_t:V_t\rightarrow[-\epsilon,\epsilon]$ be a strictly
increasing function with $\beta_t(s)=s+t$ on
$[-t-\epsilon,-t-\epsilon/2]$ and $\beta_t(s)=s-t$ on
$[t+\epsilon/2,t+\epsilon]$. Define a smooth embedding
$f_t:V_t\times N\hookrightarrow M$ by:

$$f_t(s,m)=(\beta_t(s), m).$$

Let $\bar J_t$ be the ${\partial\over \partial s}-$invariant almost
complex structure on $V_t\times N$ such that $\bar
J_t({\partial\over \partial s})=R$ and $\bar J_t|_{\xi}=J|_{\xi}$.
Glue the almost complex manifold $(M\backslash f_t(V_t\times N), J)$
to $(V_t\times N, \bar J_t)$ via $f_t$ to obtain the family of
almost complex structures $J_t$ on $M$.

 Notice that  each $J_t$ agrees with $J$ away from the collar $(-\epsilon,\epsilon)\times H$.
 And on this collar,  it agrees with $J$ on $\xi$.  Suppose $N$ is
 separating, denote $M\backslash N=W\cup U$, where $W$ has a concave
 boundary and $U$ a convex boundary.
When $i\to \infty$ the neck-stretching process results in an almost
complex structure $J_{\infty}$ on the union of symplectic
completions $\overline{W}=(-\infty,0]\times N\cup W$ of $W$ and
$\overline U=U\cup [0,+\infty)$ of $U$.  On the cylindrical ends, we
require $J|_\infty(\partial_s)=R$ and $J|_\infty=J|_\xi$ similar to
the definition of $\overline J_t$.  In an exact same way, we define
$J_\infty$ on $SN=((-\infty,+\infty)\times N, d(e^t\alpha))$, the
symplectization of $N$.

Let $M_\infty=\overline{W}\cup SN\cup\overline{U}$, and $J_\infty$
be the almost complex structure defined above.
 Let $\Sigma$ be a Riemann surface with nodes.  A \textit{level-$k$ holomorphic building} consists of the following data:
\begin{enumerate}[(i)]
 \item (level) A labelling of the components of $\Sigma\backslash\{\text{nodes}\}$ by integers $\{1,\cdots,k\}$ which are the \textit{levels}.
Two components sharing a node differ at most by $1$ in levels. Let
$\Sigma_r$ be the union of the components of
$\Sigma\backslash\{\text{nodes}\}$ with label $r$.

 \item (asymptotic matching) Finite energy holomorphic curves $v_1:\Sigma_1\rightarrow U$,
$v_r:\Sigma_r\rightarrow SN$, $2\leq r\leq k-1$, $v_k:\Sigma_k\rightarrow W$.
Any node shared by $\Sigma_l$ and $\Sigma_{l+1}$ for $1\leq l\leq
k-1$ is a positive puncture for $v_l$ and a negative puncture for
$v_{l+1}$ asymptotic to the same
 Reeb orbit $\gamma$.  $v_l$ should also extend continuously across each node within $\Sigma_l$.
\end{enumerate}

Now for a given stretching family $\{J_{t_i}\}$ as previously
described, as well as $J_{t_i}$-curves
$u_i:S\rightarrow(M,J_{t_i})$, we define the Gromov-Hofer
convergence as follows:

A sequence of $J_{t_i}$-curves $u_i:S\rightarrow(M,J_{t_i})$ is said
to be \textit{convergent to a level-$k$ holomorphic building} $v$
 in Gromov-Hofer's sense, using the above notations, if
there is a sequence of maps $\phi_i:S\rightarrow \Sigma$, and for
each $i$, there is a sequence of $k-2$ real numbers $t_i^r$,
$r=2,\cdots,k-1$, such that:

\enone[(i)]
 \item (domain) $\phi_i$ are locally biholomorphic except that they may collapse circles in $S$ to nodes of $\Sigma$,
 \item (map) the sequences $u_i\circ\phi_i^{-1}:\Sigma_1\rightarrow U$, $u_i\circ\phi_i^{-1}+t_i^r:\Sigma_r\rightarrow SH$, $2\leq r\leq k-1$, and
$u_i\circ\phi_i^{-1}:\Sigma_k\rightarrow W$ converge in
$C^\infty$-topology to corresponding maps $v_r$ on compact sets of
$\Sigma_r$. \entwo

\nono Now the celebrated compactness result in SFT reads:

\thmone[\cite{compactness}] \label{theorem:cpt} If $u_i$ has a fixed
homology class, there is a subsequence $t_{i_m}$ of $t_i$ such that
$u_{t_{i_m}}$ converges to a level-$k$ holomorphic building in the
Gromov-Hofer's sense. \thmtwo

Definitions and statements above holds true for bordered stable maps
with no extra complications, as long as the Lagrangian boundary does
not intersect the contact type boundary $N$.  Since the choice of
almost complex structure will play an important role in subsequent
sections, we would like to specify a special class of adjusted
almost complex structures for later applications.

Denote $e_1'$, $e_2'$ and $e_3'$ as the pre-images of the three
edges of $\Phi_{\CC P^2}$, numbering in a coherent way as in $F_4$
in Section \ref{subsection:description of classes}.

\defone We say $J\in\JJ^{\epsilon}_{tadj}$, the space of \textit{compatible
toric adjusted almost complex structures}, if $J$ is compactible
with $\omega_{std}|_{X_0}$, and adjusted to the hypersurface
$M_\epsilon=\Phi_{\CC P^2}^{-1}(\{x_2=1-\epsilon\})$, while $e'_1$,
$e_2'$ and $e'_3$ are $J$-holomorphic in $X_0=\Phi_{\CC
P^2}^{-1}(\{x_2<1-\epsilon\})$. Moreover, $J$ is invariant under the
circle action generated by Reeb flow in a neighborhood of
$M_\epsilon$.
\deftwo

It is not hard to see that $\tadj$ is non-empty.  Notice $e_1'$,
$e_3'$ intersects $M_\epsilon$ transversely, and are foliated by
simple orbits of the circle action.  Moreover, the Liourville vector
field near $M_\epsilon$ is invariant under the circle action, and is
tangent to $e_1'$ and $e_3'$ in a neighborhood of $M_\epsilon$.
Therefore, one only needs to define the almost complex structure to
be adjusted, whose restriction to the contact distribution is
invariant under the circle action, then extend to the rest of $X_0$
in an $\omega_{std}|_{X_0}$-compatible way so that $e_i'\cap X_0$
are holomorphic for $i=1,2,3$.

We would like to point out that, one can still achieve
transversality within $\tadj$ because no (punctured) holomorphic
curves lies entirely in the region we fixed the almost complex
structure, with the exceptions of $e_i'$, which are clearly regular
on their own right (see Wendl's automatic transversality in Section
\ref{subsubsection:D4}). Moreover, the space of such almost complex
structures is contractible, because it is just the space of sections
of a bundle with contractible fibers with prescribed values on a
closed set.

\subsection{Contributions of holomorphic disks of Maslov index $2$}

In this section we will compute terms involved in \eqref{eq:PO} by
studying evaluation of several moduli spaces. We first study the configurations of limits under
neck-stretching of holomorphic disks of Maslov index $2$, then study
all possible cases of resulting holomorphic buildings.

Here we fix some more notation convenient for our exposition. For a
Borel-Moore class $B$, we consider the moduli space of holomorphic
disks punctured at an interior point, and with one marked point on
the boundary, which we denote as $\MM^k_1(B; M, J)$ if the interior
puncture is asymptotic to $k$ times of a simple Reeb orbit. We also
consider the evaluation maps:


$$ev^i: \MM_1^{k}(B; M, J)\rightarrow \cN$$

\nono where $\cN$ is the Morse-Bott manifold where the interior
puncture lies in. When no confusion is likely to occur, we sometimes
suppress $M$ and $J$.

\subsubsection{Neck-stretching of holomorphic disks}

Given $J\in\tadj$, we may perform neck-stretching described in
\ref{section:neck-stretching}, and denote $J^+:= J_{\infty}|_{X_0}$,
 $J^-:= J_{\infty}|_{X_1}$.
Recall that $X_0$ can be compactified to $F_4$ by collapsing the
circle action on the boundary.  Under this operation, the asymptotic
boundary of $X_0$ collapses to the edge $e_4$, and (part of) $e_i'$
gives rise to $e_i$ in $F_4$ for $i=1,2,3$. Given a
$J^+$-holomorphic punctured curve $C$ of finite energy with boundary
on $L$ (possibly empty), from the asymptotic analysis in \cite{B},
$C$ can also be compactified to a well-defined $2$-cycle $\overline
C\in C_2(F_4, L;\ZZ)$ with $\overline C\cap e_4\geq0$. Also
$\overline C\cap e_i\geq 0$, $i=1,2,3$, following the positivity of
intersection and the definition of $\tadj$. Notice that $[D_i]$,
$i=1,2,3,4$ also forms a basis of $H_2(F_4,L;\ZZ)$ from
(\ref{eq:rel1})(\ref{eq:rel2}), and by Poincare duality, elements in
$H_2(F_4,L;\ZZ)$ can be identified by their pairings with divisors
$e_i$ for $i=1,2,3,4$. We therefore proved:

\lemmaone\label{lemma:decomposition} For $J\in\tadj$, an irreducible
$J^+$-holomorphic curve $C$ with finite energy, possibly with boundary on $L$
and punctures on $\partial^\infty(X_0)$, has its
class in the positive span of $\{[D_i]\}_{i=1}^4$. In particular,
the Maslov index $\mu^{\Phi}(C)\geq0$, and equality holds if and
only if $[C]=k[D_4]$ for some $k\in\ZZ^{\geq0}$. \lemmatwo

\nono We are now ready to prove:

\lemmaone\label{lemma:four classes} For $J\in\tadj$, the equation
(\ref{eq:PO}) has at most four terms of
contributions coming from $[D_1]$, $[D_2]$, $[D_3]$ and
$[D_4']$.\lemmatwo

\pfone Given $J\in\tadj$, by neck-stretching we obtain a family of
almost complex structure $J_t$.  Given a homology class $A$ which
admits $J_{t_i}$-holomorphic disks with Maslov index $2$ for a
sequence $t_i\nearrow\infty$, $i\in\ZZ^+$.  By the compactness
theorem \ref{theorem:cpt}, it converges to a holomorphic building.
We then have one of the
following cases:\\

\nono\textit{Case $1$: the $X_1$-part of the holomorphic building is
empty.}

Since $X_0$ is symplectomorphic to $F_4\backslash e_4$,
(\ref{eq:first cherns 1}) and Lemma \ref{lemma:decomposition}
implies $D_1$,
$D_2$, $D_3$ are the only possibilities, otherwise the Maslov index must exceed $2$.\\

\nono\textit{Case $2$: the $X_1$-part of the holomorphic building is
non-empty.}

Consider the $X_1$-part of the holomorphic building $S_1$.  Since it
must have periodic orbits as asymptotes, it is a Borel-Moore cycle
of class $kH'$ for some $k\in \ZZ^+$.  Therefore,
$c_1^{\Phi}(S_1)\geq 1$ by (\ref{eq:first cherns 3}),
 and the equality holds only when $k=1$.  To close up
this cycle in $\CC P^2$, one must cap $S_1$ by some cycle in $X_0$.
However, from our computations in Section \ref{subsection:CZ indices
and c1} and Lemma \ref{lemma:decomposition} we saw that all
classes but multiples of $[D_4]$
have positive first Chern number. Therefore, the only
$J_\infty$-holomorphic building with Maslov index $2$ consists of a
cycle in class $H'$ in $X_1$ and a holomorphic disk in the class
$2[D_4]$ in $X_0$. The class they form in $H_2(\CC P^2,L;\ZZ)$ is
$[D'_4]$.

\pftwo

Lemma \ref{lemma:four classes} narrows our study down to four
classes.  Notice the above two lemmata assumes no genericity of $J$.
Moreover, from the proof we see that to understand the contributions
of $[D_i]$, $i=1,2,3$, it suffices to study the stretching limit. To
understand holomorphic disks in $[D_4']$, we need a slightly more
detailed description of the limit holomorphic building:

\lemmaone\label{lemma:configuration of D4} When
$t\rightarrow\infty$, $J_t$-holomorphic disks in class $D_4'$
converge to a holomorphic building consisting of the following levels,
if existed:

\enone[(1)]
\item the $X_1$-part is a holomorphic plane in class $H'$ with one asymptotic puncture of multiplicity $2$;
\item the symplectization part is a trivial cylinder with one asymptotic puncture of multiplicity $2$ on both
positive and negative sides;
\item the $X_0$-part is a holomorphic disk in class $2[D_4]$ with a single puncture of multiplicity $2$.
\entwo \lemmatwo

\pfone In the proof of Lemma \ref{lemma:four classes} we already saw
that $X_1$-part can only be of class $H'$ and $X_0$-part is a cycle
in class $2[D_4]$ by counting Maslov indices and numbers of
punctures.  To see that $X_0$-part is a holomorphic disk with a
single puncture of multiplicity $2$ instead of two simple punctures,
notice otherwise the holomorphic building will be forced to have at
least genus $1$, since simple orbits cannot be capped by disks on
$X_1$ side. This verifies $(3)$.

In the symplectization part, since all orbits have the same period,
and the positive end has exactly one orbit of multiplicity $2$, the
negative end also has at most $2$ orbits counting multiplicities.
Again since simple orbits do not close up in $X_1$, there must be
two negative ends counting multiplicity.  Since the $\lambda$-energy
is now zero for the symplectization part, the image of the
symplectization part is a trivial cylinder.  Since branched covers
over the trivial cylinder always create genus in this holomorphic
building, we conclude that the symplectization part is indeed an
unbranched double cover of the trivial cylinder. This verifies (2),
as well as that $X_1$-part has exactly one puncture of multiplicity
$2$.  The rest of assertions in $(1)$ is easy.

\pftwo

\subsubsection{Contribution of $[D_i]$}  In this section, we prove
that:

\propone\label{prop:contribution of Di} For generic $J\in\tadj$,
$deg(ev_{0*}[\MM_1([D_i];J)])=1$ for $i=1,2,3$.\proptwo

\pfone Let us perform a
neck-stretch on $J$, so that all disks of $\MM_1([D_i];J_t)$ lie
entirely in $X_0$.  Since $J$ is cylindrical near $M_\epsilon$, we
have the following claim:

\lemmaone\label{lemma:acs lemma} $(X_0, J^+=J_\infty|_{X_0})$ is
biholomorphic to an open set $U$ of a closed symplectic manifold
with a compatible almost complex structure $(\widetilde X_0,
\widetilde \omega, \widetilde J)$, so that the following holds:

\enone[(i)]
\item $(\widetilde{X}_0,\widetilde \omega)$ is in fact the result of the symplectic cut constructed in
Section \ref{subsection:description of CP2}, i.e. $\wt X_0=X_0''=F_4$;
\item $\Sigma=\widetilde X_0\backslash U$ is a $\widetilde
J$-divisor.

\entwo \lemmatwo

This is simply a translation between the set-up of relative
invariants of \cite{LR} and the one of symplectic field theory in
the case when Reeb orbits foliates the contact type hypersurface.
$\widetilde X_0$ as a symplectic manifold comes from collapsing the
circle action on $\partial X_0$, which forms a symplectic divisor
$\Sigma$. For some small $\delta>0$, near $\Sigma$ the symplectic form of $\widetilde{X}_0$ can
be written as:

\eqone\label{eq:local
symplectic}\omega=\pi^*\tau_0+d(r'\lambda'),\eqtwo

\nono for $\delta>r'>0$. Here $\tau_0$ is a symplectic form on
$\Sigma$, $r'$ a radial coordinate; $\pi$ is the radial projection
to $\Sigma$, and $\lambda'$ a connection $1$-form (in our case it is
also a contact form) on level sets of $r'$, satisfying
$d\lambda'=\pi^*\tau_0$. Given any complex structure $J$ on
$\Sigma$, $J$ can be lifted to the horizontal distributions $\xi$
(i.e. the contact distributions), while the almost complex structure
on the whole neighborhood can be defined by further requiring
$J(r'\partial_{r'})=R'$. Here $R'$ is the Hamiltonian flow generated
by the local (in our case also global) $S^1$-action.  Conversely,
given an almost complex structure satisfying $J(r'\partial_{r'})=R'$
and invariant under the circle action on $\widetilde U\backslash
\Sigma$ where $\widetilde U$ is a neigbhorhood of $\Sigma$, it has a
natural extension to $\Sigma$.


On the SFT side, endow a symplectic form written as $d(r\lambda)$ to
the collar of $N=\partial X_0$, $1+\delta\geq r>1$, where $\lambda$
is the contact form on $N$. This coordinate can be transformed back
to the one in Section \ref{section:neck-stretching} by taking a
$log$-function on the cylindrical coordinate. The zero level-set
there becomes the level set $r=1$ in the current coordinate. In the current
coordinate, the toric adjustedness of $J^+$ is equivalent to the
invariance under both flows of $r\partial_r$ and $R$, and that
$J^+(\partial_r')=R$, where $R$ are the contact distribution and the
Reeb flow, respectively.

Notice the fact that $(N\times (1,1+\delta), d(r\lambda))$ is
symplectomorphic to $(N\times (0,\delta), d\lambda+d(r\lambda))$
just by shifting the $r$-coordinate.  By choosing $\tau_0$ so that
$\pi^*\tau_0=d\lambda$, the symplectic cut, in
perspective of this coordinate change, is simply to glue a divisor
$\Sigma$ to $(N\times (0,\delta), d\lambda+d(r\lambda))$, then the
symplectic form extends natually.  In
particular, the shift above provides a symplectic identification of a
collar neighborhood of $\partial X_0\subset X_0$ and $\Sigma$ of $\widetilde
U\backslash\Sigma$. Under such an identification, $J^+$ induces an
almost complex structure $\widetilde J$ on $\widetilde U\backslash
\Sigma$, which is invariant under $r'\partial_{r'}$ and the
Hamiltonian flow $R'$ by the assumption of toric adjustedness.  It
is then straightforward to see that $\wt J$ extends to the cut
divisor $\Sigma$ in the new coordinate.  Extend further the identification
on $\widetilde U$ to a diffeomorphism between $U=\widetilde
X_0\backslash \Sigma$ and $X_0$, we induce $\widetilde J$ by $J^+$
on the whole $\widetilde X_0$.

Given Lemma \ref{lemma:acs lemma} and the removable singularity
theorem, we may identify $\MM_1([D_i],J_t)$ with the moduli space of
holomorphic disks without punctures in $F_4$ endowed with an toric
adjusted almost complex structure $\widetilde J$ so that equivariant
divisors as $e_i$ are $\wt J$-holomorphic. A problem arises after
the compactification: $\widetilde J$ is never generic, in the sense
that $e_4$ has negative chern number, yet always $\widetilde
J$-holomorphic.  We cannot use Cho-Oh's classifications either,
because $\wt J$ is not toric. However, we can still prove:

\lemmaone\label{lemma:no redundancy Di} $\MM_1([D_i];\widetilde X_0,
\widetilde J)$ is compact for $i=1,2,3$.  Hence when $t$ is large enough,

$$ev_{0*}([\MM_1([D_i];\CC P^2, J_t)])=ev_{0*}([\MM_1([D_i];\widetilde X_0, \widetilde J)]), i=1,2,3.$$

\lemmatwo

\pfone To understand the left hand side, we may consider the problem
in the limit and replace the left hand side by $X_0$ and $J^+$. By
the same analysis as in Lemma \ref{lemma:four classes}, case 2, any
leveled curves coming as a limit for $t\rightarrow+\infty$ must have
empty $X_1$-part. Further Lemma \ref{lemma:acs lemma} identifies
such a curve as one on the right hand side which has no
$e_4$-component.  Hence the conclusion follows provided one can prove
$[D_i]$ are indecomposable on the right hand side. Corresponding classes are
also indecomposable on the left hand side by a similar reasoning.  But in our application
we will assume $L$ is monotone so all classes on the left hand side
 are clearly indecomposable for $\w$-area reason,
so we will omit the actual proof.


Take $[D_1]$ as an example, and the rest of the cases are similar.
Assume $u:\Sigma\rightarrow(\widetilde X_0, \widetilde J)$ is a
stable curve in the moduli space of right hand side with irreducible
components $\Sigma_1,\dots,\Sigma_k$, and the homology classes of
all of these components are written in terms of basis
$\{[D_i]\}^4_{i=1}$.  Let $\Sigma_1\dots,\Sigma_l$ be components of
$e_4$,  while $[e_4]=[D_1]+[D_3]-4[D_4]$.  Now Lemma
\ref{lemma:decomposition} implies that $[\Sigma_j]$, $j>l$ all lie
in the positive cone of spanned by $[D_i]$ for $i=1,2,3,4$.  By
comparing coefficients of $[D_3]$, we conclude that $l=0$.  It then
follows easily that $k=1$ and $[\Sigma_1]=[D_1]$. Since $[D_1]$
pairs trivially with $[e_4]$, our claim is proved by positivity of
intersections.

\pftwo

Lemma \ref{lemma:no redundancy Di} implies that
$\MM_1([D_i];\widetilde X_0, \widetilde J)$ is in fact compact even
with no genericity assumption since the class itself is
indecomposable. Therefore, the standard cobordism arguments apply.
In particular, one may choose a generic path $\{J_t\}_{t\in[0,1]}$
connecting $J_0$ and $J_1=\widetilde J$ for $J_0$ also satisfying
that $e_i$, $i=1,2,3,4$ are $J_0$-holomorphic. Recall from \cite{CO}
that there is an integrable complex structure $J_0$ where
$ev_{0*}([\MM_1([D_i];\widetilde X_0, J_0)])$ is known to be $[L]$,
hence concluding our proof of Proposition \ref{prop:contribution of
Di}.

\pftwo

\subsubsection{The contribution of $[D'_4]$}\label{subsubsection:D4}

Our goal of this section is to prove:

\propone\label{prop:contribution of D4} For generic choice of $J\in
\tadj$,
$$deg(ev_{0*}[\overline \MM_1(D_4';\CC P^2, J)])=2.$$
\proptwo

As already explained in previous sections, we only need to
consider $J\in \tadj$ and its neck stretched sufficiently long.  We
first briefly review Wendl's automatic transversality theorem.

One of the new ingredients of Wendl's theorem is the introduction of
the invariant \textit{parity}, defined in \cite{propertyII}, to the
formula. Let $Y$ be a symplectic cobordism, where $Y^{\pm}$ are the
positive (resp. negative) boundaries. Given a $T$-periodic orbit
$\gamma$ of $Y^{\pm}$, one has an associated asymptotic operator,
which takes the form of $\mathbb{A}= -I_0\partial_t-S(t)$ on
$L^2(S^1,\RR^2)$ by taking a trivialization of the normal bundle.
Here $I_0$ is the standard complex structure on $\RR^2$, while
$S(t)$ is a continuous family of symmetric matrices.  For
$\lambda\in \sigma(\mathbb{A})$, one may define a \textit{winding
number} $w(\lambda)$ to be the winding number of nontrivial
$\lambda$-eigenfunction of $\mathbb{A}$.  It is proved in
\cite{propertyII} that $w(\lambda)$ is an increasing function of
$\lambda$ which takes every integer value exactly twice.  For
 non-degenerate operators $\mathbb{A}$ (i.e. $0\notin \sigma(\mathbb{A})$), we define

$$\alpha_+(\mathbb{A})=max\{w(\lambda)|\lambda\in\sigma(\mathbb{A}), \lambda<0\},$$
$$\alpha_-(\mathbb{A})=min\{w(\lambda)|\lambda\in\sigma(\mathbb{A}), \lambda>0\}.$$

\nono and the parity $p(\mathbb{A})=\alpha_+(\AAA)-\alpha_-(\AAA)
(\text{mod }2)$.  If $\AAA$ is degenerate, we define
$\alpha_\pm(\AAA\pm\delta)$ and $p(\AAA\pm\delta)$ for small
$\delta>0$.  For a given puncture, the actual perturbation depends
on which of $Y^\pm$ it lies on, as well as whether the moduli space
we consider constrains the puncture inside a Morse-Bott family.
Chris Wendl pointed out to the author that, in our case when the
contact type boundary is foliated by a $2$ dimensional family of
Reeb orbits, since the eigenvalue $0$ has multiplicity $2$,
\textit{either way of perturbation incurs odd parity}.

Now given a non-constant punctured holomorphic curve
$u:\Sigma_g\rightarrow Y$, the virtual index is computed as:

$$ind(u)=(n-3)\chi(\Sigma)+2c_1(u)+\sum_{\gamma^+}(\mu_{CZ}(\gamma^+)+\frac{1}{2}dim(\cN))-
\sum_{\gamma^-}(\mu_{CZ}(\gamma^-)-\frac{1}{2}dim(\cN)),$$

\nono Here $\gamma^\pm$ runs over all positive (resp. negative)
punctures, and $\cN$ is the Morse-Bott manifold formed by the Reeb
orbits.  We now define the \textit{normal chern number} as:

$$2c_N(u)=ind(u)-2+2g+\#\Gamma_0+\#\pi_0(\partial\Sigma_g).$$

\nono Here, $\Gamma_0$ denotes the number of punctures of even
parities, hence in our applications when the contact type boundary is foliated
by $2$ dimensional family of Reeb orbits, this term
always vanishes.  Having understood these, Wendl's automatic transversality
theorem reads:

\thmone[\cite{WenAT}]\label{thm:Wendl} Suppose \textup{dim}$Y$ = 4
and $u:(\Sigma, j)\rightarrow (Y,J)$ is a non-constant curve with
only Morse-Bott asymptotic orbits. If
$$ind(u) > c_N(u) + Z(du),$$
\nono then u is regular.
\thmtwo

For the contribution to \eqref{eq:PO} of holomorphic disks in class
$D'_4$, we consider the gluing problem of $\MM^2(H';X_1,J^-)$
and $\MM^2_1(2[D_4]; X_0, J^+)$.  Note that this is sufficient
by the configuration analysis of the limit holomorphic building in
Lemma \ref{lemma:configuration of D4}. The standard gluing argument
requires the following conditions:

\enone[(1)]
\item Curves in both $\MM^2(H';X_1,J^-)$ and $\MM^2_1(2[D_4]; X_0, J^+)$ are regular;
\item $ev^1\times ev^1:\MM^2(H';X_1,J^-)\times \MM^2_1(2[D_4]; X_0, J^+)\rightarrow S^2\times S^2$
is transversal to the diagonal $\Delta\subset S^2\times S^2$.  Here
$S^2=\cN$ is exactly the Morse-Bott family parametrizing Reeb orbits on
$M_\epsilon$. \entwo

One sees that item $(2)$ is automatic since the first component of
the evaluation map is surjective onto $S^2$. This corresponds to the
standard fact in Gromov-Witten theory that, given any compatible
almost complex structure $J$ in $\CP^2$, an embedded $J$-conic
$\Sigma$ and a point $p\in\Sigma$, there is a unique $J$-complex line tangent to
$\Sigma$ at $p$. For $(1)$ we apply
Wendl's automatic transversality in dimension $4$.

The virtual index of an irreducible curve $C\in\MM^2(2[D_4]; X_0,
J^+)$ reads:

$$ind(u)=(2-3)(2-1-1)+0+0-(0-1)=1.$$

The computation also shows that, for generic $J$, the
compactification of this moduli space does not contain irreducible
curves with critical points or sphere bubbles since these are
codimension $2$ phenomena. On the other hand, we can compute
$c_N(u)=0$. Therefore, automatic transversality holds for all
 $C\in \MM^2_1(2[D_4]; X_0, J^+)$.

To show that disk bubbles do not appear, we again use Lemma
\ref{lemma:acs lemma} to identify the moduli space
$\MM^2_1(2[D_4];X_0,J^+)$ to one on $F_4$, denoted $\MMc$,
where $\widetilde J$ is the extended almost complex structure.

\defone
$\MMc$ is the moduli space of $\widetilde J$-holomorphic disks
$u:(D^2,j)\rightarrow (\widetilde X_0, \widetilde{J})$ which
satisfies the following:

\begin{itemize}
\item $u$ has an interior marked point $x$ and a boundary marked point $y$,
\item $u(\partial D)\subset L$, $u(x)\in e_4$, $du(x)=0$ with order $1$.
\end{itemize}
\deftwo

Now by collapsing the Reeb orbits on $\partial X_0$, a stable
punctured disk in $\MM^2_1(2[D_4];X_0,J^+)$ descends to a
stable disk in $\MMc$.  The order of vanishing of $du$ exactly
corresponds to the multiplicity of the asymptotic Reeb orbit.

\lemmaone\label{lemma:D4 transversality} Holomorphic disks
$u\in\MM_1^{\text{crit}}(2[D_4];e_4,\widetilde X_0,L;\widetilde J)$
are regular for generic $\widetilde J$. Moreover, the moduli space
is compact. \lemmatwo

\pfone The argument is almost word-by-word taken from the case of
open manifolds. Since $u$ cannot develop more critical points other
than $x$ for generic choice of $\widetilde J$, we may apply Wendl's
automatic transversality, Theorem \ref{thm:Wendl}. We have the
Fredholm index:

$$\text{ind}(u)=-1+2\cdot 2=3,$$
$$c_N(u)=\frac{1}{2}(3-2+1)=1.$$

Since we have a unique critical point of order $2$,
$$3=\text{ind}(u)>c_N(u)+Z(du)=1+1=2,$$

\nono verifying the transversality of $u\in\MMc$.  We now only need
to show that $\partial\MM_1^{\text{crit}}(2[D_4];e_4,\widetilde
X_0,L;\widetilde J)=\emptyset$.  The argument of Lemma \ref{lemma:no
redundancy Di} shows that the only possible type of reducible stable
curve $u\in\partial\MM_1^{\text{crit}}(2[D_4];e_4,\widetilde
X_0,L;\widetilde J)$ consists of a union of $2$ disks in class
$[D_4]$ (by comparing coefficients of either $[D_1]$ or $[D_3]$ for
possible irreducible decompositions). However, given a sequence
$u_k\in\MM_1^{\text{crit}}(2[D_4];e_4,\widetilde X_0,L;\widetilde
J)$ converging to $u$, $u_k(x)\in e_4$ are always critical values
which cannot approach the boundary.
 If a disk bubble occurs, one of the components inherits such a critical
 point, thus has intersection index with $e_4$ at least $2$.  But
 this contradicts the fact that each component is in class $[D_4]$.

\pftwo

To compute the evaluation map, we now may choose a generic path
connecting $\{J^s\}_{s\in[0,1]}$ connecting $J^1=\widetilde J$ and
the standard toric complex structure $J^0$ of $F_4$ as in \cite{CO},
while requiring $e_i$, $i=1,2,3,4$ are $J^s$-holomorphic. In view of
the arguments in Lemma \ref{lemma:D4 transversality}, the moduli
space $\MM_1^{crit}(2[D_4]; F_4, L; J^s)$ does not develop disk or
sphere bubbles. Moreover, $2[D_4]$ does not admit a multiple cover
more than $2$-fold, whereas these $2$-fold covers are in fact curves
in $\MMc$ instead of the boundary of the moduli space
$\partial\MM_1^{\text{crit}}(2[D_4];e_4,\widetilde X_0,L;J^s)$.
Therefore, the standard cobordism argument in \cite{MSJ} applies. In
particular, $ev_{0*}[\MMc]=ev_{0*}[\MM_1^{crit}(2[D_4]; F_4, L;
J^0)]$. By Cho-Oh's classification, $\MM_1^{crit}(2[D_4]; F_4, L;
J^0)$ consists only of double covers of embedded disks in
$\MM_1([D_4]; F_4, L; J^0)$ with critical points on intersections
with $e_4$. Therefore, from the identification of $\MM^2_1(2[D_4];X_0,J^+)$ and $\MMc$,

$$\begin{aligned}
ev_{0*}[\MM^2_1(2[D_4];X_0,J^+)]&=ev_{0*}[\MMc]\\
                                &=ev_{0*}[\MM_1^{crit}(2[D_4]; F_4, L; J^0)]=2[L].\end{aligned}$$

On the $X_1$ side, what concerns us is $\MM^2(H';X_1,J^-)$. We
already saw from the argument of item $(2)$ that these curves
one-one correspond to closed
 curves in $\CC P^2$ of line class which are tangent to the given embedded conic.  In particular no bubbling
 or critical points occurs for these curves.  The virtual index of such a curve $C_1$ is:

$$ind(C_1)=(2-3)(2-1)+2+1-0=2,$$

\nono and $c_N(C_1)=0$.  This verifies item $(1)$.  Moreover,
$ev^1:\MM^2(H;X_1,J^-)\rightarrow S^2$ is surjective and of
degree $1$ from considering the Gromov-Witten invariants of tangent
lines on the conic after closing up the orbits on the boundary.
Therefore, the standard gluing arguments apply and leads to the
following commutative diagram:

$$\xymatrix{\MM^2(H';X_1,J^-) \,\,{}_{ev^1}\times_{ev^1}
\MM^2_1(2[D_4]; X_0, J^+)
\ar[rr]^(.6){glue}\ar[dr]_{\widetilde{ev}_0} & & \mathcal M_1(D'_4; \CC P^2, L) \ar[dl]^{ev_0} \\
&L&}$$

\nono Here $\widetilde{ev}_0$ is the evaluation of $\MM^2_1(2[D_4];
X_0, J^+)$ to the boundary marked points.  It then follows
that, for $t$ sufficiently large,

\eqone\label{eq:mapping degree of D4'} deg[ev_0:\MM_1(D'_4; \CC P^2,
L; J_t)\rightarrow L]=2. \eqtwo

\section{Completion of the proof}\label{s:completion}

To summarize, we have computed evaluation maps of all holomorphic
disks of Maslov index $2$ of $(\CC P^2, J_t)$ when $t$ is
sufficiently large. Plugging these inputs into(\ref{eq:PO}) we
deduce that the potential function of $\CC P^2$ in the toric
degeneration picture is written as:

$$\PO^u(y)=T^{u_1}y_1+T^{u_2}y_2+T^{4-u_1-4u_2}y_1^{-1}y_2^{-4}+2T^{2-2u_2}y_2^{-2}.$$

\nono By requiring $u_1=u_2=\frac{2}{3}$, the equation of critical
points are written as:

\begin{equation} \label{eq:1}
\left\{ \begin{aligned}
         & y_1^2y_2^4=1 \\
         & 1-4y_2^{-5}y_1^{-1}-4y_2^{-3}=0
                          \end{aligned} \right.
                          \end{equation}

\nono We therefore deduce that $y_1y_2^2=\pm 1$.  When we take it to
be $1$, we have $y_2^3=8$ thus clearly has three solutions. This
verifies that $\mathbf{u}=(\frac{2}{3}, \frac{2}{3})$ is indeed a
critical point for some values of $y_1$ and $y_2$. Moreover, $QH_*(\CC P^2,\Lambda)$ is semi-simple and
decomposes into $3$ direct factors of $\Lambda$. Writing $QH^*(\CC
P^2;\Lambda)=\Lambda[z]/(z^3-T)$, the idempotents are simply
${\displaystyle
\frac{1}{3}\epsilon_i^{-2}(z^2+\epsilon_iz+\epsilon_i^2)}$ for
$i=1,2,3$.  Here $\epsilon_i$ are the roots of $x^3-T=0$ in
$\Lambda$.  Now Corollary \ref{cor:main lemma} implies
$L(\mathbf{u})$ is indeed superheavy with respect to some symplectic
quasi-state.  This concludes our proof to Theorem \ref{thm:main}.

\rmkone Our example manifests two interesting aspects of the effect
of the choice of Novikov rings.  Namely, we found $3$ local systems
on the exotic monotone fiber, which is different from the case of
the calculation of \cite{FOOOdeg}, where the monotone exotic fiber
in $S^2\times S^2$ only has half number of local systems of the
standard monotone fiber, i. e. the product of equators.

According to comments due to Kenji Fukaya, combining results from
\cite{AFOOO}, our computation implies that this single exotic fiber
is sufficient to generate certain Fukaya category with
characteristic zero coefficients.  However, this fiber is disjoint
from $\RR P^2$, thus cannot generate any version of Fukaya category
with characteristic $2$ coefficients (in fact our torus is always a
zero object for characteristic $2$ Fukaya categories of $\CC P^2$).
This shows that the choice of the characteristic of coefficient rings could be more than
technical. \rmktwo

\brmk Leonid Polterovich brought up another very interesting
question:  is it possible to distinguish the symplectic
quasi-states/morphisms for the three idempotents of
$QH^*(\CP^2;\Lambda)$?  These three quasi-states/morphisms are
intuitively very closely related as is pointed out also by
Polterovich.  Given the Novikov field
$\Lambda_{nov}=\{a_iT^{\lambda_i}: \lambda_i\in \Z,
\lim\lambda_i=+\infty\}$, then $\Lambda_{nov}[z]/(z^3-T)$ is already
a field. However, when we tensor this ring by $\Lambda$, the
algebraic closure of $\Lambda_{nov}$, the ring is only semi-simple,
as indicated in our computation.  Therefore, the identity in the
field $\Lambda_{nov}[z]/(z^3-T)$ splits into $3$ idempotents after a
purely algebraic procedure, thus intuitively, the three symplectic
quasi-states/morphisms are ``algebraically splitted" from the
original quasi-state/morphism as well.

It was known to Entov and Polterovich \cite{EP1} that for $S^2$
symplectic quasi-states are unique up to certain normalization.
However, the corresponding statement for quasi-morphism is not know
even for the spectral quasi-morphisms in this case.  For $\CP^n$,
$n\geq2$, there is no results available.  It would be very
interesting to understand how essential are these algebraic
extensions of symplectic quasi-states/morphisms.

\ermk

\vskip 2mm

\nono Department of Mathematics,  Michigan State University\\
East Lansing, MI48824, USA\\

\nono \textsf{wwwu@math.msu.edu}

\end{document}